\documentclass[11pt, reqno]{amsart}
\usepackage{amsmath, amsthm, amscd, amsfonts, amssymb, graphicx,tikz, color,xcolor,mathtools}
\usepackage[all]{xy}
\usepackage{hyperref}
\usepackage[T1]{fontenc} 
\usepackage{mathrsfs}
\usepackage{ wasysym}
\usepackage[inline]{enumitem}
\usepackage{hyperref}
\usepackage{multicol}
\usepackage{rotating}
\usepackage{graphicx,fancyhdr}
\usepackage{multirow}

\hypersetup{colorlinks=true,citecolor=cyan,linkcolor=red,linktocpage=true}

\newtheorem{theorem}{Theorem}[section]
\newtheorem{lemma}[theorem]{Lemma}
\newtheorem{coro}[theorem]{Corollary}

\newtheorem{prop}[theorem]{Proposition}

\theoremstyle{definition}
\newtheorem{remark}[theorem]{Remark}

\newtheorem{nota}[theorem]{Notations}
\newtheorem{definition}[theorem]{Definition}

\newtheorem{stepp}{Claim}
\newtheorem{stepo}{Step}

\newtheorem{claim}{Claim}
\newtheorem{claimm}{Claim}
\newtheorem{case}{Case}
\newtheorem{subcase}{Subcase}[case]

\newcommand{\pf}{\begin{proof}}
\newcommand{\epf}{\end{proof}}

\newcommand{\yd}[1]{{}^{ #1 }_{ #1 }\mathcal{YD}}

\newcommand{\Der}{\operatorname{Der}}
\newcommand{\Derinn}{\operatorname{InnDer}}

\newcommand{\co}{\operatorname{co}}

\newcommand{\N}{{\mathbb N}}

\newcommand{\ku}{{\Bbbk}} 
\newcommand{\kku}{{\ku[z^p]}}

\newcommand{\toba}{{\mathscr B}} 

\newcommand{\Z}{{\mathbb Z}} 
\newcommand\id{\operatorname{id}} 
 
\newcommand{\reg}{S}

\newcommand{\ch}{{\rm char\,}}
\newcommand{\Pc}{\mathcal P}

\newcommand{\B}{{\mathcal B}}
\newcommand{\Zc}{\mathcal Z}

\newcommand{\Aut}{\operatorname{Aut}} 
\newcommand{\Auts}{\Aut_{\operatorname{super}}}
\newcommand{\Auth}{\Aut_{\operatorname{Hopf}}} 
\newcommand{\Res}{\operatorname{Res}}
\newcommand{\imm}{\operatorname{im}}

\newcommand{\GK}{\operatorname{GKdim}}

\newcommand{\J}{\mathbf{J}}

\begin{document}
\title[On the super Jordan plane]{On the super Jordan plane}
\author[Andruskiewitsch]{Nicol\'as Andruskiewitsch}
\address[N.~Andruskiewitsch]{CIEM-CONICET. 
Medina Allende s/n (5000) Ciudad Universitaria, C\'ordoba, Argentina.
\newline Department of Mathematics and Data
Science, Vrije Universiteit Brussel, Pleinlaan 2, 1050 Brussels, Belgium}
\email{nicolas.andruskiewitsch@unc.edu.ar}
\author[Dumas]{Fran\c cois Dumas}
\address[F.~Dumas]{Universit\'e Clermont Auvergne - CNRS, LMBP - UMR 6620. \newline
F-63000 Clermont-Ferrand, France}
\email{Francois.Dumas@uca.fr}

\begin{abstract}
We study algebraic properties of the super Jordan plane~$\B$ introduced by I. Angiono, I. Heckenberger and the first named author. Concretely we show that $\B$ is super-prime and has a
super-simple super-artinian ring of fractions. We also compute the groups of superalgebra and braided Hopf algebra  automorphisms. 
Instrumental for our approach is the embedding of the Jordan plane as a subalgebra of $\B$.
\end{abstract}

\thanks{
\newline
N.~A. was  partially supported by the Secyt (UNC),
CONICET (PIP 11220200102916CO) and FONCyT-ANPCyT (PICT-2019-03660. He  is grateful  to 
 the LMBP of the Universit\'e Clermont Auvergne for the hospitality 
 during his visit in June 2025.}

\date {September 1, 2025}
\maketitle
\setcounter{tocdepth}{1}
\tableofcontents

\section*{Introduction}
\subsection*{The Jordan plane}
This is a quadratic algebra generated by $\mathbf{x}$
and $\mathbf{y}$ with defining relation
 \begin{equation*}
\mathbf{y}\mathbf{x}-\mathbf{x}\mathbf{y}=-\frac12\mathbf{x}^2.
\end{equation*}
It was considered in the literature in 
relation with AS-regular algebras \cite{ArS},
 quantum groups \cite{G,DMMZ, zak,LM, ohn} and
 the classification of pointed Hopf algebras with finite 
 Gelfand-Kirillov  dimension ($\GK$) \cite{AAH1,AAH2};
 indeed it is a Nichols algebra in characteristic $0$. 
 The name \emph{Jordan} is due to \cite{DMMZ}.
 Algebraic properties of the Jordan plane were studied
 in \cite{AlD,BLO,Iyudu,Shirikov0,Sh2007}, notice that it is a domain of $\GK$  $2$ and a Koszul algebra.

\subsection*{The super Jordan plane}
This is the graded algebra $\B$ generated by  $x,y$ (in degree 1)
with two defining relations 
\begin{align*}
x^2 &= 0, &  yz-zy-xz &= 0,
\end{align*}
where $z=xy+yx$. Thus, one relation is quadratic and one is cubic.
In characteristic $0$,  $\B$  is indeed a Nichols algebra. 
 It was introduced in \cite{AAH2}, being essential in the
 classification of Nichols algebras of finite $\GK$ over abelian groups. It has $\GK$  $2$ but it is not a domain.
 Further properties of $\B$ were established:
 
 \medbreak
 \noindent $\circ$ The  classification of the simple finite-dimensional, respectively indecomposable of dimension 2  and 3, $\B$-modules was given in \cite{ABDF17},
 assuming $\ch  \Bbbk = 0$.
 
  \medbreak
 \noindent $\circ$ The  classification of the simple finite-dimensional modules over a  suitable bosonization of the  super Jordan plane 
 was given in \cite{ABDF},
  assuming $\ch \Bbbk = 0$.
 
  \medbreak
 \noindent $\circ$ 
 The Hochschild homology and cohomology of $\B$ were studied in \cite{RS} assuming $\ch\ku=0$ and $\ku$ algebraically closed.
 So, in particular, the Lie algebra of derivations can be described
 from \emph{loc. cit.}
 
 \medbreak
 \noindent $\circ$ 
The graded braided commutativity of the Hochschild cohomology 
for a class of braided Hopf algebras, including $\B$, 
was proved in \cite{CS}.

 \medbreak
 \noindent $\circ$ 
 The Drinfeld double of a suitable bosonization of $\B$
was introduced in \cite{APP}; it gives rise naturally to 
a Hopf superalgebra which is an extension of a super commutative one by the enveloping  algebra of $\mathfrak{osp}(1|2)$.

\subsection*{Main results}
In the present paper we provide further information on
$\B$. Namely, we show that 
 $\B$ is a superalgebra, meaning that it has a $C_2$-grading $\B = \B_0 \oplus \B_1$. 
 Strictly speaking, one should speak of $C_2$-graded algebras, coalgebras, 
 etc., and reserve the term super to structures requiring the super 
 transposition  like Lie algebras, Hopf algebras and so on. But 
 at the end of the day both settings are related, so we speak in a
 uniform manner of super algebras, super ideals, etc.  
Our main results are:

\medbreak
\noindent $\bullet$ 
The algebra $\B$ has infinite global dimension (Corollary \ref{cor:gldim});
a different proof follows from \cite{RS}.

\medbreak
\noindent $\bullet$ 
The zero divisors of $\B_0$
and the homogeneous zero divisors of $\B$ are determined
(Propositions \ref{evenZD} and \ref{oddZD}).

\medbreak
\noindent $\bullet$
 The algebra $\B_0$ can be localized at
the set $\reg _0$ of its regular elements  
and  $Q_0\coloneqq \B_0\reg _0^{-1}$ is semisimple artinian (Theorem \ref{QB0}).

\medbreak
\noindent $\bullet$
The superalgebra $\B$ is a super-prime ring (Theorem
\ref{Bsp}). It has not been apparent to us whether the algebra $\B$ is prime
so we let this as an open question.

\medbreak
\noindent $\bullet$
The algebra $\B$ can be localized at the set $\reg $ of homogeneous regular elements in $\B$ 
and the ring of fractions $Q\coloneqq \B\reg ^{-1}$ has
a $C_2$-grading whose even part is $Q_0=\B_0\reg _0^{-1}$. 
We have $Q=\B\reg _0^{-1}$. Moreover
the superalgebra $Q$ is a super-simple super-artinian ring
(Theorem \ref{QB}).   

\medbreak
\noindent $\bullet$ The group $\Auts\B$ of automorphisms 
of the superalgebra $\B$, i.e., the algebra automorphisms
that preserve the $C_2$-grading, is computed (Theorem \ref{AutB}).
It turns out that $\Auts \B \simeq N \rtimes \ku^{\times}$,
where the automorphisms in $N$ are exponentials of locally nilpotent derivations (Proposition \ref{prop:locnilp}).

\medbreak
\noindent $\bullet$ The group of automorphisms 
of the braided Hopf algebra $\B$, with the structure introduced  in \cite[Section~3.3]{AAH2}, see Section \ref{section:braided-Hopf}, is computed assuming characteristic 0 (Proposition \ref{prop:hopf-automorphisms}).

\subsection*{Notations and conventions}
The (commutative) base field is denoted by~$\ku$.
We set $ p \coloneqq \ch\ku$ which is assumed $\not=2$; additional constraints  will be imposed accordingly.
As usual in group theory, $C_2 = \{\overline{0}, \overline{1}\}$ is the group of order 2. When there is no danger of confusion, we simply set $C_2 = \{0, 1\}$.

Let $R$ be an algebra. The center of $R$ is denoted by $Z(R)$;
if $R$ is a superalgebra, then the supercenter of $R$
is denoted by $\Zc(R)$.
If $\varphi$ is an automorphism of $R$, then $R^{\varphi}$
denotes the subalgebra of points fixed by $\varphi$.

\section{The super Jordan plane as an iterated Ore extension}
\begin{definition}\label{defsJp}The super Jordan plane is the algebra $\B$ presented by generators $x,y$ with defining relations
\begin{equation}\label{relB}
x^2=0,\quad yz-zy-xz=0, \ \text{with notation} \ z=xy+yx.
\end{equation}
It follows from these relations that 
\begin{equation}\label{relBbis}
xz=zx.
\end{equation}
The algebra $\B$ is $\Z$-graded: 
$\B = \oplus_{n \geq 0} \B_n$, where $x$ and $y$
have degree $1$.
Hence $\B$ is a $C_2$-graded algebra (or a superalgebra) by 
\begin{align*}
\B_{\overline{0}} &= \oplus_{n \text{ even}\ } \B_n, &
\B_{\overline{1}} &= \oplus_{n \text{ odd}\ } \B_n.
\end{align*}
\end{definition}

\begin{nota} We consider the commutative algebra
\begin{equation}A\coloneqq \ku[x]/(x^2),\end{equation} 
known as the algebra  of dual numbers,
and the  polynomial  algebra $A[z]$. Any element $f=\sum_{i=0}^N(\alpha_i+\beta_i x)z^i \in A[z]$, where
$\alpha_i, \beta_i \in \ku$, can be uniquely written as
\begin{align}\label{f+-}
f  =& f^++xf^-  \text{ with}& f^+&=\sum_{i=0}^N\alpha_i z^i \in\ku[z],
& f^- &=\sum_{i=0}^N\beta_i z^i \in\ku[z].
\end{align}

\medbreak
The Euler derivative in the polynomial ring $\ku[z]$ 
is denoted by $D$; namely
\begin{align}\label{Euler}
D(f)&\coloneqq z\partial_z(f) & \text{for any } f &\in\ku[z].
\end{align}
Then $\ker D= \kku$, where of course  $\kku = \ku$ when $p=0$.
\end{nota}

The algebra $\B$ can be presented as an Ore extension.

\begin{prop}\label{Bore}
\begin{enumerate}[leftmargin=*,label=\rm{(\roman*)}]
\item \label{item:Bore}
There exist an automorphism $\tau$ and a $\tau$-derivation $d$ 
of the algebra $A[z]$ such that the super Jordan plane is isomorphic 
to the Ore extension
\begin{equation}\label{BOre} A[z][y\,;\,\tau,d].
\end{equation}
The automorphism $\tau$ and the $\tau$-derivation $d$ are defined by
\begin{align}\label{taudxz}
\tau(x)&=-x, & \tau(z) &=z, & d(x)&=z,  & d(z)&=xz.
\end{align}

 \item \label{item:basis}
 The set $\mathcal X = \{x^az^by^c : a\in\{0,1\}, b,c\in\N_0\}$ is a $\ku$-basis of $\B$. 
\end{enumerate}
\end{prop}

Part  \ref{item:basis} was claimed without proof in  \cite[Lemma  3.5]{APP}.

\pf We define $\tau:A\to A$ by $\tau(a)=\alpha-\beta x$ for any $a=\alpha+\beta x$, $\alpha,\beta\in\ku$. An easy calculation (using $x^2=0$) proves that $\tau(ab)=\tau(a)\tau(b)$ for all $a,b\in A$. Clearly $\tau$ defines an automorphism of $A$, which can be  canonically extended to $A[z]$ setting $\tau(z)=z$. Using the decomposition \eqref{f+-} we have \begin{equation}\label{tau+-}
\tau(f)=f^+-xf^-\ \text{for any } f=f^++xf^-\in A[z].
\end{equation}

We define $d:A\to A$ in two steps. First, we set
\begin{equation}\label{dEuler}
d(f)=xD(f)=D(f)x\ \text{for any } f\in\ku[z].
\end{equation}
 
Then using \eqref{f+-} we define
\begin{equation}\label{d+-}
d(f)=xD(f^+)+zf^-\ \text{for any } f=f^++xf^-\in A[z].
\end{equation}
If $g=g^++xg^- \in A[z]$, we have,  again because $x^2=0$, that
\begin{align*}
d(fg)
&=d(f^+g^++x(f^+g^-+f^-g^+))\\
&=xD(f^+g^+)+z(f^+g^-+f^-g^+)\\
&=xD(f^+)g^++xf^+D(g^+)+zf^+g^-+zf^-g^+\\
&=(f^+-xf^-)(xD(g^+)+zg^-)+(xD(f^+)+zf^-)(g^++xg^-)\\
&=\tau(f)d(g)+d(f)g.
\end{align*}
Hence $d$ is a $\tau$-derivation of $A[z]$.
Now, there  exists an algebra map 
$\psi:\B \to A[z][y,\tau, d]$ sending the generators $x$ and $y$ to their homonimous. It was proved in \cite[Proposition 3.5]{AAH2}
that  $\mathcal X$ linearly generates $\B$.
Clearly $\psi(\mathcal X)$ is a basis of $A[z][y,\tau, d]$, 
hence $\psi$ is an isomorphism and $\mathcal X$ is a basis of $\B$.
\epf

\begin{coro} \label{cor:gldim}
The algebra $\B$ is Noetherian and has infinite global dimension.
\end{coro}

\pf
The first claim follows from \cite[Theorem 2.6]{GW}.
By Proposition \ref{Bore} \ref{item:basis} the set $\mathcal X=\{x^az^by^c : a\in\{0,1\}, b,c\in\N_0\}$ is a $\ku$-basis of $\B$. 
Now $\mathcal X=\mathcal X_0\cup\mathcal X_1$ where 
$\mathcal X_0=\{1,x\}$ and $\mathcal X_1=\{x^az^by^c : a\in\{0,1\}, b,c\in\N_0, (b,c)\not=(0,0)\}$, hence
$\B=A\oplus C$ where  $C$ is the $\ku$-span of $\mathcal X_1$.
Clearly $C$ is a left $A$-module.
From the commutation relations in $\B$ we see that $C$ is also a right $A$-module.
Indeed, it is evident that $z^bx=xz^b$ hence $z^bx\in C$ for any $b\geq 1$.
By straightforward calculations, see
\cite[Lemma 3.6]{APP}, we have $y^cx\in C$ for any $c\geq 1$.
Thus, finally,  $z^by^cx\in C$ for any $z^by^c\in \mathcal X_1$.
We deduce that $A$ is an $A$-bimodule direct summand of $\B$. 
Since $\B$ is Noetherian and free over $A$, we conclude by 
\cite[Theorem 7.2.8]{mcconnell-robson}.
\epf

\begin{remark} 
Since $A$ is not a domain, the degree related to $y$ in the Ore extension \eqref{BOre} does not satisfy the usual property on the degree of a product. For instance we have in $\B$  the relation
\begin{equation}\label{relxyzn}(xy)^{n+1}=z^nxy\quad\text{for any $n\in\N_0$},
\end{equation}
which follows by induction from $(xy)(xy)=x(-xy+z)y=xzy=zxy$ and implies in particular that
\begin{equation}\label{relxyz}(z-xy)xy=xy(z-xy)=0
\end{equation}
Zero divisors in $\B$ are discussed in Propositions \ref{evenZD} and \ref{oddZD}.
\end{remark}

\begin{lemma}\label{lemdtau}We have in $A[z]$ the following properties:
\begin{enumerate}[leftmargin=*,label=\rm{(\roman*)}]
\item
 $\tau^2=\id$ and $\tau d=-d\tau$,

\medbreak
\item
 $\ker d=\kku$, and $A[z]^\tau=\ku[z]$,

\medbreak
\item $d^2$ is a derivation.
\end{enumerate}
\end{lemma}
\pf Clearly $\tau^2=\id$, by \eqref{tau+-}. By \eqref{f+-}, \eqref{tau+-}, \eqref{dEuler} and \eqref{d+-} we have 
\begin{align*}
&\begin{aligned}
d\tau(f)&=d(f^+-xf^-)=d(f^+)-\tau(x)d(f^-)-d(x)f^-\\
&=d(f^+)+x(xD(f))-zf^-=xD(f^+)-zf^-,\\
\tau d(f)&=\tau(xD(f^+)+zf^-)=-x\tau(D(f^+))+z\tau(f^-)\\
&=-xD(f^+)+zf^-,
\end{aligned}& \text{for any } f &\in A[z],
\end{align*}which proves (i). 
By \eqref{d+-} we have  
$d(f)=0$ if and only if $D(f^+)=0$ and $f^-=0$, hence by \eqref{Euler} if and only if $f=f^+\in\kku$.
Using \eqref{tau+-} we have
\begin{equation}f-\tau(f)=2xf^-
\end{equation}then $\tau(f)=f$ if and only if $f=f^+\in\ku[z]$ and point (ii) is proved. Here we used that $\ch\ku\not=2$.

\medbreak
Let $f,g$ be two elements of $A[z]$. We have $d(fg)=\tau(f)d(g)+d(f)g$, then
\begin{align*}
d^2(fg)&=
\tau^2(f)d^2(g)+(d\tau+\tau d)(f)d(g)+d^2(f)g\\
&=fd^2(g)+d^2(f)g \ \ \ \text{by point (i)}
\end{align*}
which shows that $d^2$ is an ordinary derivation of $A[z]$.
\epf

\section{The Jordan plane as a subalgebra}
Following the notations of \cite{AAH1,AAH2} we recall that the Jordan plane $\J$ is the algebra presented by generators $\mathbf{x}$ and $\mathbf{y}$
with defining relation
 \begin{equation}\label{oldnotJ}
 \mathbf{y}\mathbf{x}-\mathbf{x}\mathbf{y}=-\frac12\mathbf{x}^2.
 \end{equation}

\begin{prop}\label{B+}
The formal differential operator algebra $A[z][y^2\,;\,d^2]$ is a subalgebra of $\B$, and its subalgebra $\ku[z][y^2\,;\,d^2]$
can be identified with the Jordan plane $\J$ via the isomorphism given by
\begin{equation}\label{eq:jordan-subalgebra}
\mathbf{x} \mapsto z \quad \text{and} \quad \mathbf{y} \mapsto-\frac12y^2.
\end{equation}
\end{prop}

\pf We know that $yf=\tau(f)y+d(f)$ for any $f\in A[z]$. Hence we have
\begin{equation}\label{eq:comm-rules}
y^2f=\tau^2(f)y^2+(d\tau+\tau d)(f)y+d^2(f)=fy^2+d^2(f),
\end{equation}
by Lemma \ref{lemdtau} (i); this proves the first assertion. By \eqref{taudxz} we have 
\begin{equation*}d^2(z)=d(xz)=-xd(z)+d(x)z=-x^2z+z^2=z^2,\end{equation*}
then the subalgebra $\ku[z][y^2\,;\,d^2]$ is generated over $\ku$ by $z$ and $y^2$ with relation 
$y^2z-zy^2=z^2$ and the proof is complete.\epf

\begin{remark}\label{rksJ}
 By iteration of the commutation rules \eqref{eq:comm-rules} we have 
\begin{equation}\label{commJ}
y^{2n}f=\sum_{i=0}^n\binom{n}{i}d^{2n-2i}(f)y^{2i}\qquad\text{for any } f\in\ku[z].
\end{equation}
\end{remark}

Recall that $p = \ch \ku \neq 2$, so that $\ku[z^p, y^{2p}] = \ku$ when $p=0$.

\begin{prop}\label{centerB}
 $Z(\B) = Z(\J) = \ku[z^p, y^{2p}]$.
\end{prop}

It is known that $Z(\J) = \ku[z^p, y^{2p}]$, see \cite[Theorem 2.2]{Sh2007}. 
That $Z(\B)=\ku$  for algebraically closed $\ku$, $\ch \ku = 0$, 
is proven in \cite[Theorem 3.12]{RS}.

\pf We shall prove that $Z(\B)=Z(\J)$; then \cite[Theorem 2.2]{Sh2007}
implies the statement.
Let ${b}=\sum_{i=0}^nf_iy^i$ an element of $\B$ with $f_i\in A[z]$. We have
\begin{align*}
{b}y-y{b}&=\sum_{i=0}^nf_iy^{i+1}-\sum_{i=0}^n(\tau(f_i)y+d(f_i))y^i\\
&=d(f_0)+\sum_{j=1}^{n}(f_{j-1}-\tau(f_{j-1})-d(f_j))y^j+(f_n-\tau(f_n))y^{n+1}.
\end{align*}
Then ${b}$ commutes with $y$ if and only if 
\begin{align*}d(f_0)&=0,\\
d(f_j)&=f_{j-1}-\tau(f_{j-1})\ \ \text{for }1\leq j\leq n, \\
f_n-\tau(f_n)&=0.
\end{align*}
From Lemma \ref{lemdtau} (ii) 
we deduce that $f_0\in\kku$. Then $f_0-\tau(f_0)=0$, thus $f_1\in\kku$ and  inductively $f_j\in\kku$ for all $0\leq j\leq n$. We deduce that the centralizer of $y$ in $\B$ is $\kku[y]$.\smallskip

Take ${b}=\sum_{i=0}^nf_iy^i$ in this centralizer, with $f_i\in \kku$, $f_n\not=0$ and suppose that 
${b}\notin\kku$.
We have
\begin{equation*}
{b}x-x{b}=(f_ny^n+\cdots)x-x(f_ny^n+\cdots)
=f_n(\tau^n(x)-x)y^n + \rho
\end{equation*}
where  $\deg_y \rho < n$.
If ${b}$ commutes with $x$, then $\tau^n(x)=x$, hence $n=2\ell$ is even. Let $m = \sup \{i \text{ odd} : 1\leq i\leq n-1 \text{ and } f_i\not=0\}$.
Then
\begin{equation*}
{b}x-x{b}=\sum_{i=\frac{m+1}2}^{\ell}f_{2i}(y^{2i}x-xy^{2i})+f_{m}(y^{m}x-xy^{m})+ \varrho
\end{equation*}
where   $\deg_y \varrho <m$.
The first term of this sum belongs to $A[z][y^2\,;\,d^2]$ by Proposition \ref{B+}, while the second summand is of the form
$-2f_mxy^{m}+\cdots$ since $m$ is odd. If ${b}x=x{b}$, then $f_m=0$,
a contradiction. We deduce that  if ${b}\in\kku[y]$ commutes with $x$, then ${b}\in\kku[y^2]$. In this  case, ${b}$ commutes with $z$ if and only if ${b} \in Z(\J)$. 
\epf

\medbreak
 We use the notations $\ku[z][y^{2p}]=\ku[z]$ and $\ku[z^p][y^2]=\ku[y^2]$ when $p=0$.

\begin{lemma}\label{defetaD} \emph{(i)} There exists an automorphism $\eta$ of $\J$ satisfying 
\begin{align}\label{eq:defeta}
\eta(y^2) & =y^2-z & &\text{and}& \J^\eta &=~\ku[z][y^{2p}].
\end{align}
 Hence we  have  the following equalities:
\begin{align}\label{eta}
xf &=\eta(f)x, & zf &=\eta(f)z,&  xyf &=\eta(f)xy & \text{ for any }f &\in\J.
\end{align}

\emph{(ii)} There exists a linear endomorphism $\nabla $ of $\J$ satisfying 
\begin{align}\label{eq:defDelta}
\ker \nabla  & =\kku[y^2]  & &\text{ and } &\nabla (z) &=z.
\end{align}

 Hence we have  the following equalities:
\begin{align}\label{Delta}yf&=fy+\nabla (f)x,& 
y^2f &=fy^2+\nabla (f)z,&  \text{for any } f&\in\J.
\end{align}
\end{lemma}

\pf It is well known, see Remark \ref{compAutBJ} below,
that the automorphisms of $\J$ are of the form 
$z\mapsto \alpha z$ and $y^2\mapsto\alpha y^2+p(z)$ 
for some $\alpha\in\ku^\times$ and $p(z)\in\ku[z]$. 
Let $\eta$ be the automorphism of $\J$ corresponding 
to $\alpha=1$ and $p(z)=-z$. Then relations \eqref{eta} 
translate commutations rules $xz=zx$, $zy^2=y^2z-z^2$, 
$(xy)z=z(xy)$ and $(xy)y^2=(y^2-z)xy$ deduced 
from \eqref{relB} and \eqref{relBbis}. 
We have for $\eta$ the following generic expression 
\begin{equation}\label{defeta}
\eta(f)=\sum_{i=0}^na_i(z)(y^2-z)^{i} \quad \text{for any } f=\sum_{i=0}^na_i(z)y^{2i}\in\J.
\end{equation}
By induction we  immediately have $(y^2-z)^i=y^{2i}-izy^{2(i-1)}$. Hence
\begin{equation}\label{defetabis}
\eta(f)-f=-z\sum_{i=1}^nia_i(z)y^{2(i-1)} \quad \text{for any } f=\sum_{i=0}^na_i(z)y^{2i}\in\J.
\end{equation}
We deduce that $\eta(f)=f$ if and only if $k\in\ku[z]$ when $\ch\ku=0$, and $f\in\ku[z,y^{2p}]$ when $\ch\ku=p$. The proof of (i) is complete.\smallskip

Let us define the linear map $\nabla :\J\to\J$ by
\begin{equation}\label{defDelta}
\nabla (f)=\sum_{i=0}^nD(a_i(z))(y^2-z)^{i} \quad \text{for any } f=\sum_{i=0}^na_i(z)y^{2i}\in\J.\end{equation}
Then we have
\begin{align*}yf&=\sum_{i=1}^n(a_iy+d(a_i))y^{2i}=fy+\sum_{i=1}^nD(a_i)xy^{2i}
\\& = fy+\sum_{i=1}^nD(a_i)(y^2-z)^{i}x =fy+\nabla (f)x,
\end{align*}
Since $\ker D=\kku$ by \eqref{Euler}, it is clear by \eqref{defDelta} that $\ker\nabla =\kku[y^2]$.\epf

\begin{remark}\label{rkDelta}Let us observe that the linear map $\nabla $ is not a derivation of~$\J$. More precisely we have in $\J$ 
\begin{align*}
\nabla (fg)z&=y^2(fg)-(fg)y^2=(y^2f-fy^2)g+f(y^2g-gy^2)\\
&=\nabla (f)zg+f\nabla (g)z=\nabla (f)\eta(g)z+f\nabla (g)z;
\end{align*}
that is, 
$\nabla (fg)=f\nabla (g)+\nabla (f)\eta(g)$ for all $f,g\in\J$. The restriction of $\nabla $ to~$\ku[z]$ is just the Euler derivative $D=d^2$ of \eqref{Euler}.
\end{remark}

\section{The super Jordan plane as a superalgebra}
Keep the notations above. Recall that $\B$ is $C_2$-graded.
Below we set
$\B_0 = \B_{\overline{0}}$ and $\B_1= \B_{\overline{1}}$.
The identification \eqref{eq:jordan-subalgebra} of the Jordan plane with a subalgebra of $\B$ appeared already in \cite {ABDF17}, where it was observed that 
$1,x,y,xy$ generate the $\J$-module $\B$, see \cite[Proposition 2.3]{ABDF17}. More precisely, we have:

\begin{prop}\label{Bsuper}The $C_2$-grading
of the  superalgebra $\B$ is determined by assigning parity $0$ to the elements of $\J$ and parity 1 to $x$ and to $y$, that is
\begin{align}\label{B0B1}
\B &=\B_0\oplus\B_1& &\text{where}& 
\B_0 &=\J\oplus\J xy && \text{and} & \B_1 &=\J x\oplus\J y.
\end{align}
\end{prop}
\pf The linear decomposition 
$\B=\J\oplus\J y\oplus \J x  \oplus\J xy$ follows from Proposition \ref{B+} and the relation $x\J=\J x$, which is deduced from \eqref{eta}.
Then  \eqref{B0B1} follows at once.
\epf

\begin{remark}
The first part of the  previous proof shows that $\B$ is a 
free left $\J$-module of rank 4.
\end{remark}

\begin{remark}\label{rem:algebra-B_0}
The algebra structure of $\B_0 = \J\oplus\J xy$ is determined by the algebra structure of $\J$ and the identities:
\begin{align*}
xyf &= \eta (f)xy, \quad f \in \J,&&\text{and}& (xy)^2 &= zxy.
\end{align*}
\end{remark}

The following lemma gives an expression of the product of two homogeneous elements.

\begin{lemma}\label{twohomo}Let $f_0=f_0^++f_0^-xy$ and $g_0=g_0^++g_0^-xy$ be two even elements, and $f_1=f_1^+x+f_1^-y$ and $g_1=g_1^+x+g_1^-y$ two odd elements, with $f_0^\pm,g_0^\pm,f_1^\pm,g_1^\pm\in\J$. Then we have
\begin{align}
\label{twoeven}
f_0g_0&=f_0^+g_0^++\Bigl(f_0^+g_0^-+f_0^-\eta(g_0^++g_0^-z)\Bigr)xy\\
\label{twoodd}
f_1g_1&=f_1^-(g_1^+z+g_1^-y^2)+\Bigl(f_1^+\eta(g_1^-)+f_1^-\bigl(\nabla (g_1^-)-g_1^+\bigr)\Bigr)xy\\
\label{evenodd}
f_0g_1&=\Bigl(f_0^+g_1^++f_0^-\eta(g_1^+z+g_1^-y^2)\Bigr)x+(f_0^+g_1^-)y\\
\label{oddeven}
g_1f_0&=\Bigl(g_1^+\eta(f_0^+)+g_1^-\big(\nabla (f_0^+)-f_0^-\eta(y^2)\bigr)\Bigr)x+g_1^-(f_0^++f_0^-z)y
\end{align}
\end{lemma}
\pf This follows from the equalities
\begin{align*}
xyx&=zx,& (xy)^2 &=zxy,& xy^2&=(y^2-z)x& \text{and } yxy &=(z-xy)y
\end{align*}
that are deduced easily from the relations 
\eqref{relB}, \eqref{eta} and \eqref{Delta}. 
\epf

It is well-known  that  $\J$ is a domain; however,
 $\B_0$ is not one, see \eqref{relxyz}.

\begin{prop}\label{B0sp}The subalgebra $\B_0$ is a semiprime, 
but not prime, ring.
\end{prop}

\pf Let $f_0=f_0^++f_0^-xy\in\B_0$, where $f_0^{\pm} \in \J$,
such that $f_0\B_0f_0=0$; in particular $f_0^2=0$. Relation \eqref{twoeven} implies 
\begin{align*}
(f_0^+)^2 &=0& &\text{and}& f_0^+f_0^-+f_0^-\eta(f_0^++f_0^-z)=0.
\end{align*}
in the domain $\J$. Then $f_0^+=0$ hence $f_0^-\eta(f_0^-)z=0$ which implies $f_0^-=0$. Hence $f_0=0$ which proves the first assertion, see \cite[Theorem 3.7]{GW}.

To prove the second assertion, we choose $f_0=xy$ and $g_0=z-xy$. 
For any $h_0=h_0^++h_0^-xy$ with $h_0^+,h_0^-\in\J$, we have
\begin{align*}
f_0h_0g_0=xy(h_0^++h_0^-xy)(z-xy)=xyh_0^+(z-xy)=\eta(h_0^+)xy(z-xy)=0
\end{align*}
 using \eqref{relxyz}. Hence $f_0\B_0g_0=0$ and $\B_0$ is not prime.
\epf

The supercenter $\Zc(R) = \Zc_0(R) \oplus\Zc_1(R)$
of a superalgebra $R = R_0 \oplus R_1$ is 
the super-subspace of elements that super-commute with
those in $R$, i.e.,
\begin{align*}
\Zc_0(R)  &= R_0 \cap  Z(R),\\
\Zc_1(R)  &= \{f \in R_1: fg = (-1)^{j}gf \text{ for all } g \in R_j, j\in C_2\}.
\end{align*}

Once again, we recall that $\ku[z^p, y^{2p}] = \ku$ when $p = \ch \ku = 0$.

\begin{prop}\label{supercenterB}  $\Zc(\B)  = \ku[z^p, y^{2p}]$.
\end{prop}

\pf
Let $f=f_0+f_1 \in \Zc(\B)$, with $f_0\in\Zc_0(\B) 
= \B_0 \cap Z(\B)  = \ku[z^p, y^{2p}]$, by Proposition \ref{centerB}. Write $f_1=f_1^+x+f_1^-y\in\Zc_1(\B)$
where $f_1^{\pm} \in\J$. As $x^2=0$, the
relation $f_1x=-xf_1$ reduces to  
$ f_1^-yx=-xf_1^-y$  which is equivalent, by \eqref{eta}, to
\begin{align*}
&f_1^-(-xy+z)  =-\eta(f_1^-)xy, 
\\ \intertext{ which implies }
 &f_1^-z +(\eta(f_1^-)-f_1^-)xy =0,
\end{align*}
hence $f_1^-=0$. Then the equality $f_1y=-yf_1$ becomes
\begin{align*}
f_1^+xy=-yf_1^+x=-(f_1^+y+\nabla (f_1^+)x)x=-f_1^+yx=f_1^+xy-f_1^+z,
\end{align*}
thus $f_1^+=0$,  $f_1=0$ and $f=f_0\in\ku[z^p, y^{2p}]$.
\epf


\section{The superalgebra of fractions}
The goal of this subsection is to localize the algebras $\B_0$ and $\B$
at the regular homogeneous elements, i.e., those that are not zero divisors.
So we start by computing these.

\begin{prop}\label{evenZD}Let $f_0 \in\B_0$. 
The following assertions are equivalent.

\begin{enumerate}[leftmargin=*,label=\rm{(\roman*)}]
\item\label{item:evenZD1} $f_0$ is a right zero divisor in $\B$,

\medbreak
\item\label{item:evenZD2} $f_0$ is a left zero divisor in $\B$,

\medbreak
\item\label{item:evenZD3} $f_0$ is a right zero divisor in $\B_0$,

\medbreak
\item\label{item:evenZD4}$f_0$ is a  left zero divisor in $\B_0$,

\medbreak
\item\label{item:evenZD5} $f_0$ belongs to the union of the  ideals $\J xy=xy\J$ and $\J (z-xy)=(z-xy)\J$  of $\B_0$.
\end{enumerate}\end{prop}

\pf  The implications (iii) $\Rightarrow$ (i) and (iv) $\Rightarrow$ (ii) are evident.
We start by the following claim.

\begin{claim}\label{claim1}
The statements (iii), (iv) and (v)  are equivalent.
\end{claim} 
It follows from \eqref{eta} that $\J xy=xy\J$ and $\J (z-xy)=(z-xy)\J$; since  $(xy)^2=zxy$, these  are  ideals of $\B_0$. By \eqref{relxyz} we get the implications (v) $\Rightarrow$ (iii) and (v) $\Rightarrow$ (iv).

Suppose that  $f_0=f_0^++f_0^-xy$ is a right zero divisor in~$\B_0$ where 
$f_0^{\pm} \in~\J$;
we may assume that $f_0 \neq 0$. 
There exist $g_0^{\pm} \in \J$, not simultaneously~$0$,  such that $(f_0^++f_0^-xy)(g_0^++g_0^-xy)=0$. 
By \eqref{twoeven} we have  the  relations
\begin{equation}\label{zerodiv}
f_0^+g_0^+=0\quad\text{and}\quad f_0^+g_0^-+f_0^-\eta(g_0^+)+f_0^-\eta(g_0^-)z=0
\end{equation}
in the domain~$\J$. 
If $f_0^+=0$, then $f_0 = f_0^-xy \in \J xy$.
If $f_0^+\not=0$, then $g_0^+=0$ and $g_0^-\not=0$.  Relation \eqref{zerodiv} becomes 
\begin{align*}
0 = f_0^+g_0^-+f_0^-\eta(g_0^-)z = f_0^+g_0^-+f_0^-zg_0^-,
\end{align*}
hence $f_0^+=-f_0^-z$ and $f_0 = -f_0^-(z - xy) \in \J (z-xy)$. 
We have proved that 
\begin{align}\label{eq:zerodiv}
(f_0^-xy) \left(g_0^-(z-xy) \right) &=0&
&\text{or}& \left(f_0^-(z-xy) \right) (g_0^-xy) &=0.&
\end{align}

Thus (iii) implies (v). 
Analogous calculations exchanging the roles of $f_0$ and $g_0$  proves that that (iv) implies (v). The claim follows.

\begin{claim} \label{claim2}
(i) implies (iii).
\end{claim}
Suppose  that $f_0=f_0^++f_0^-xy$, $f_0^{\pm} \in \J$, is a right zero divisor in~$\B$. Pick $g=g_0+g_1\not=0$ with $g_0\in\B_0$ and $g_1\in\B_1$ such that $f_0g=0$. 
Then $f_0g_0=0$ in $\B_0$ and $f_0g_1=0$ in $\B_1$. If $g_0\not=0$, then $f_0$ is a zero divisor in $\B_0$. So we may assume  that $g_0=0$, that is, $g=g_1=g_1^+x+g_1^-y\not=0$ with $g_1^{\pm}\in\J$. Relation \eqref{evenodd} implies
\begin{equation*}
f_0^+g_1^-=0\qquad\text{and}\qquad f_0^+g_1^++f_0^-\eta(g_1^+z+g_1^-y^2)=0
\end{equation*}
in the domain $\J$. 
If $f_0^+=0$, then  $f_0=f_0^-xy$ is a zero divisor in~$\B_0$ by Claim~\ref{claim1}.
If $f_0^+\not=0$, then $g_1^-=0$, $g_1^+ \neq 0$ and  
\begin{align*}
0= f_0^+g_1^++f_0^-\eta(g_1^+)z = (f_0^++f_0^-z)g_1^+.
\end{align*}
Hence $f_0^+=-f_0^-z$  and $f_0 = -f_0^-(z - xy) \in \J (z-xy)$. 
The claim follows from Claim \ref{claim1}.

\begin{claim} \label{claim3}
(ii) implies (iv).
\end{claim}

The proof is analogous to the proof of Claim \ref{claim2}. \epf

\begin{prop}\label{oddZD}Let $f_1=f_1^+x+f_1^-y \in \B_1$, 
with $f_1^{\pm} \in\J$. The following assertions are equivalent.
\begin{enumerate}[leftmargin=*,label=\rm{(\roman*)}]
\item\label{item:oddZD1} $f_1$ is a right zero divisor in $\B$,

\medbreak
\item\label{item:oddZD2}$f_1$ is a left zero divisor in $\B$,

\medbreak
\item\label{item:oddZD3} $f_1$ is nilpotent,

\medbreak
\item\label{item:oddZD4} $f_1^2 = 0$,

\medbreak
\item\label{item:oddZD5}$f_1^-=0$ or 
$f_1^-y^2=-f_1^+z$.
\end{enumerate}\end{prop}

\pf 
The implications (iii) $\Rightarrow$ (i), (iii) $\Rightarrow$ (ii) 
and (iv) $\Rightarrow$ (iii) are evident.
Also, we may assume that $f_1\neq 0$.
We start by the following claim.

\begin{claimm}\label{claimm1}
(iv) and (v) are equivalent.
\end{claimm}

Applying \eqref{twoodd} for $g_1=f_1$ we obtain 
\begin{equation}\label{squareodd}
f_1^2=f_1^-(f_1^+z+f_1^- y^2)+\Bigl(f_1^+\eta(f_1^-)+f_1^-\bigl(\nabla (f_1^-)-f_1^+\bigr)\Bigr)xy.
\end{equation}
If $f_1^2=0$, then $f_1^-(f_1^+z+f_1^- y^2)=0$ in the domain $\J$ 
so the conditions of (v) are necessary.
Conversely it is clear that  $f_1^-=0$ implies  $f_1^2=0$.
Suppose then that $f_1^-\not=0$ and $f_1^-y^2=-f_1^+z$. 
We compute using \eqref{eta} and \eqref{Delta}:
\begin{align*}\Bigl(f_1^+\eta(f_1^-)+f_1^-\bigl(\nabla (f_1^-)-f_1^+\bigr)\Bigr)z
&=f_1^+\eta(f_1^-)z+f_1^-\nabla (f_1^-)z - f_1^-f_1^+z\\
&=f_1^+zf_1^- +f_1^-\nabla (f_1^-)z - f_1^-f_1^+z  \\
&=-f_1^-y^2f_1^-+ f_1^-\nabla (f_1^-)z + (f_1^-)^2y^2\\
&=f_1^-\Bigl(-y^2f_1^-+\nabla (f_1^-)z+f_1^-y^2\Bigr)
= 0.\end{align*}
Since $\J$ is a domain, it follows that $f_1^+\eta(f_1^-)+f_1^-\bigl(\nabla (f_1^-)-f_1^+\bigr)=0$. Therefore $f_1^2=0$ and the claim follows.

\begin{claimm}\label{claimm2}
If $f_1g_1=0$ for some $g_1 \in \B_1$, $g_1 \neq 0$, then $f_1^2 = g_1^2 = 0$.
\end{claimm} 
Let  $g_1=g_1^+x+g_1^-y \in \B_1 \backslash0$, $g_1^{\pm} \in \J$, 
such that $f_1g_1=0$.
If $f_1^-=0$, then $f_1^2 = 0$ by \eqref{squareodd} as we already observed. 
Also, $f_1^+\eta(g_1^-) = 0$ by \eqref{twoodd}. But $f_1^+ \neq 0$
by assumption, hence $g_1^- = 0$ and a fortiori $g_1^2 = 0$
by Claim \ref{claimm1}.

\medbreak
If $f_1^-\not=0$, then we deduce from \eqref{twoodd} that 
\begin{equation}\label{eq:claimm2}
g_1^+z+g_1^-y^2=0
\end{equation} 
hence $g_1^2=0$  by Claim \ref{claimm1}. 
Looking at the coefficient of $xy$ in \eqref{twoodd} we see that 
\[f_1^+\eta(g_1^-)-f_1^-g_1^++f_1^-\nabla (g_1^-)=0.\]
Using this identity, as well as \eqref{eta}, \eqref{Delta}
and \eqref{eq:claimm2}, we conclude that
\begin{align*}
(f_ 1^+z+f_1^-y^2)g_1^-&=f_ 1^+zg_1^-+f_1^-y^2g_1^-\\
&=f_ 1^+zg_1^-+f_1^-g_1^-y^2+f_1^-\nabla (g_1^-)z\\
&=f_ 1^+\eta(g_1^-)z-f_1^-g_1^+z+f_1^-\nabla (g_1^-)z\\
&=\bigl(f_1^+\eta(g_1^-)-f_1^-g_1^++f_1^-\nabla (g_1^-)\bigr)z=0.
\end{align*}
Since $g_1\not=0$,  $g_1^-\not=0$ by  \eqref{eq:claimm2}. 
Hence $f_1^-y^2=-f_1^+z$, thus $f_1^2 = 0$
by Claim~\ref{claimm1}.  

\begin{claimm}\label{claimm3}
If $f_1g_0=0$ for some $g_0 \in \B_0$, $g_0 \neq 0$, then $f_1^2 = 0$.
\end{claimm} 

Let $g_0=g_0^++g_0^-xy\in\B_0 \backslash 0$, $g_0^{\pm} \in \J$, 
such that $f_1g_0=0$. If $f_1^-=0$, then $f_1^2 = 0$. 

Assume that $f_1^-\not=0$. Relation \eqref{oddeven},
with the roles of $f$ and $g$ exchanged,
implies that $g_0^+=-g_0^-z$ hence $g_0^-\not=0$; 
also, it shows that
\begin{equation}\label{eq:claimm3}
f_1^+\eta(g_0^+)+f_1^-\nabla (g_0^+)-f_1^-g_0^-(y^2-z)=0.
\end{equation} 
Using Remark \ref{rkDelta} and \eqref{eq:defDelta} we get that
\[
\nabla (g_0^+)=-\nabla (g_0^-z)=-\nabla (g_0^-)\eta(z)-g_0^-\nabla (z)=-\nabla (g_0^-)z-g_0^-z.
\]
 Then the equality \eqref{eq:claimm3} becomes 
\begin{equation*}
-f_1^+\eta(g_0^-)z-f_1^-\bigl(g_0^-y^2+\nabla (g_0^-)z\bigr) =0;
\end{equation*} 
thus $f_1^+zg_0^- + f_1^-y^2g_0^-=0$. Hence
$f_1^-y^2=-f_1^+z$ and $f_1 ^2 =0$ by Claim \ref{claimm1}.

\begin{claimm}\label{claimm4}
If $g_0f_1=0$ for some $g_0 \in \B_0$, $g_0 \neq 0$, then $f_1^2 = 0$.
\end{claimm} Let $g_0=g_0^++g_0^-xy\in\B_0 \backslash 0$, $g_0^{\pm} \in \J$, 
such that $f_1g_0=0$. 
Exchanging the roles of $f$ and $g$ we have by \eqref{evenodd} the relations\begin{equation*}
g_0^+f_1^-=0\qquad\text{and}\qquad g_0^+f_1^++g_0^-\eta(f_1^+z+f_1^-y^2)=0.
\end{equation*}
If $f_1^-=0$, then $f_1^2 =0$.  If $f_1^-\not=0$, then $g_0^+=0$,  $g_0^-\not=0$ and $f_1^+z+f_1^-y^2=0$. Hence  $f_1^2 =0$  by Claim \ref{claimm1}.

\begin{claimm}\label{claimm5}
(i) implies (iv) and (ii) implies (iv).
\end{claimm} 
Let $f_1\in\B_1$ be a right (respectively left) zero divisor in~$\B$. There exists $g=g_0+g_1\not=0$ with $g_0\in\B_0$ and $g_1\in\B_1$ such that $f_1g=0$ (respectively $gf_1=0$).  Then $f_1g_1=0$ in $\B_0$ and $f_1g_0=0$ in $\B_1$ (respectively  $g_1f_1=0$ in~$\B_0$ and $g_0f_1=0$ in $\B_1$). We conclude by Claims \ref{claimm2},
\ref{claimm3} and \ref{claimm4}
that $f_1^2 = 0$.  The claim and the proposition follow.
\epf

\begin{nota}\label{notaS}We introduce the following subsets of $\B$ whose elements are explicitly described by Propositions \ref{evenZD} and \ref{oddZD}:
\begin{align*}
\Omega_0&\coloneqq \{\text{zero divisors in} \ \B_0\},\\
\reg_0&\coloneqq \B_0\setminus\Omega_0=\{\text{regular elements in} \ \B_0\},\\
\Omega&\coloneqq \{\text{homogeneous zero divisors in} \ \B\},\\
\reg &\coloneqq (\B_0\cup\B_1)\setminus\Omega=\{\text{homogeneous regular elements in} \ \B\}.
\end{align*}
Any right zero divisor being a left zero divisor in $\B_0$, it is clear that $f_0g_0\in\Omega_0$ implies $f_0\in\Omega_0$ or $g_0\in\Omega_0$ for all $f_0,g_0\in\B_0$; similarly for $\Omega$. 
It follows that
\begin{align*}
&\reg _0  \text{ is a multiplicative subset of } \B_0,
\\
&\reg   \text{ is a multiplicative subset of } \B_0\cup\B_1.
\end{align*}
\end{nota}

\begin{theorem}\label{QB0}
The set $\reg _0$ of regular elements in $\B_0$ satisfies the left and right Ore conditions in $\B_0$ and the total ring of fractions $Q_0\coloneqq \B_0\reg _0^{-1}$ is semisimple artinian.\end{theorem}
\pf 
The ring $\B_0$ is Noetherian, being a left and right $\J$-module of dimension 2 with $\J$ Noetherian.
Then  the result follows from Proposition \ref{B0sp} by direct application of \cite[Th\'eor\`eme 3.6.12]{Dix}.
\epf

\begin{remark}
The theory of homogeneous localization for graded rings has been developed by various authors and we refer to the monography \cite{NvO1} or the survey \cite{BKM}. In what follows, we rely on this general theory in the particular context of a $\Z/2\Z$-graduation, agreeing to replace the prefix {\it gr-} with the prefix {\it super-} as usually in other terminologies on superalgebras. For example, we will speak of super-prime or super-semiprime rings. Graded versions of some main structure theorems can be found in \cite[Theorem 2.10.10]{NvO1} or \cite[Theorem 6]{BKM} (Wedderburn-Artin's Theorem), 
\cite[Theorem 8.4.5]{NvO1} or \cite[Theorem 35]{BKM} (Goldie's Theorem). 
\end{remark}

\begin{remark}
Recall that the homogeneous Ore conditions are equivalent to the classical Ore conditions (\cite[Lemma 8.1.1]{NvO1} or \cite[Proposition 60]{BKM}). 
\end{remark}

\begin{remark} Recall that a superalgebra $R$ is 
left $e$-faithful if for any  $0 \neq x \in R$ homogeneous
there exists $x' \in R$ homogeneous
such that $0 \neq x'x \in R_e$. Right $e$-faithfulness is defined 
accordingly. Then $R$ is $e$-faithful if it is left and right $e$-faithful. See  \cite[Definition 7]{BKM}
or \cite[page 39]{NvO1}.

\medbreak
It is easy to see that $\B$ is $e$-faithful. It is an open question whether the algebra $\B$ is prime.
In the super sense, we have:
\end{remark}

\begin{theorem}\label{Bsp}
The superalgebra $\B$ is a super-prime ring.
\end{theorem}
\pf We have to prove that if homogeneous elements 
$f_i$ and $g_j$, where $i,j\in C_2$, satisfy $f_i\B g_j=0$, 
then $f_i=0$ or $g_j=0$, see \cite[Definition 11]{BKM} or \cite[page 64]{NvO1}. 
The proof naturally breaks down into four cases.\smallskip

\begin{case}\label{case1}
$f_0,g_0\in\B_0$ satisfy $f_0\B g_0=0$.
\end{case}  In particular, $f_0g_0=0$;  by  \eqref{eq:zerodiv} two subcases are possible.

\begin{subcase}
$f_0=f_0^+xy$ and $g_0=g_0^+(z-xy)$ with $f_0^+,g_0^+\in\J$.
\end{subcase}
Using  \eqref{eta} and the already observed equality
$xyx = zx$, we calculate 
\begin{align*}f_0xg_0&=f_0^+xyxg_0^+(z-xy)=f_0^+zxg_0^+(z-xy)\\
&=f_0^+z\eta(g_0^+)x(z-xy)=f_0^+z\eta(g_0^+)zx.\end{align*}
Then $f_0xg_0=0$ in $\B_1$ implies $f_0^+z\eta(g_0^+)z=0$ in $\J$. Hence $f_0^+=0$ or $g_0^+=0$
and finally $f_0=0$ or $g_0=0$.

\begin{subcase}
$f_0=f_0^+(z-xy)$ and $g_0=g_0^+xy$ with $f_0^+,g_0^+\in\J$.
\end{subcase}  
We compute
\begin{align*}f_0yg_0&=f_0^+(z-xy)yg_0^+xy=f_0^+(z-xy)(g_0^+y+\nabla (g_0^+)x)xy\\
&=f_0^+(z-xy)g_0^+yxy=f_0^+\eta(g_0^+)(z-xy)yxy.\end{align*}
As an auxiliary calculation, we compute using \eqref{eta}:
\begin{equation*}(z-xy)yxy=(z-xy)(-xy+z)y=(z^2-zxy)y
= z^2y -z\eta(y^2)x.
\end{equation*}
Then $f_0yg_0=0$ in $\B_1$ implies 
$f_0^+\eta(g_0^+)z^2= f_0^+\eta(g_0^+)z\eta(y^2)=0$ in $\J$ and we conclude as in the previous subcase.

\begin{case}
$f_1,g_1\in\B_1$ satisfy $f_1\B g_1=0$.
\end{case}   
We start considering the condition $f_1x g_1=0$. Writing as usual $f_1=f_1^+x+f_1^-y$
and $g_1=g_1^+x+g_1^-y$ with $f_1^{\pm}, g_1^{\pm} \in \J$, we have
\begin{align*}
0 &= f_1x g_1 = f_1^-yx (g_1^+x+g_1^-y) = f_1^-yx g_1^-y = f_1^- \eta(g_1^-) yxy
\\
&=f_1^- \eta(g_1^-) \left(-y^2x +yz\right) = f_1^- \eta(g_1^-) \left((-y^2 +z) x +zy\right),
\end{align*}
hence $f_1^-  g_1^- = 0$ in $\J$. 
\begin{subcase}
 $f_1^- = 0$. 
\end{subcase}
Here we consider the condition $f_1zg_1=0$. Since
 \begin{align*}
 0 = f_1zg_1=f_1^+xz(g_1^+x+g_1^-y)=f_1^+z\eta(g_1^-)xy,
 \end{align*}
  we have $f_1^+\eta(g_1^-)=0$ in $\J$. If $f_1^+=0$, then $f_1=0$. 
  
  If $f_1^+\not=0$ then $g_1^-=0$. Considering now the condition
$f_1yg_1=0$, we have
\begin{align*}
0 &= f_1yg_1=f_1^+xyg_1^+x=f_1^+x(g_1^+y+\nabla (g_1^+)x)x
\\
&=f_1^+\eta(g_1^+)xyx=f_1^+\eta(g_1^+)zx.
\end{align*}
We conclude that $f_1^+\eta( g_1^+)=0$ in $\J$, hence $g_1=0$.

\begin{subcase}
 $f_1^- \neq  0$, so that  $g_1^- = 0$. 
\end{subcase}
 The condition $f_1zg_1=0$ gives 
 \begin{align*}
 0 &= (f_1^+x+f_1^-y) z g_1^+x =  f_1^-yzg_1^+x = f_1^-(zy+zx)g_1^+x
 \\
 &=f_1^-zyg_1^+x=f_1^-z  g_1^+yx = f_1^-z  g_1^+z - f_1^-z  g_1^+ xy,
 \end{align*}
which implies that $g_1^+=0$ since $f_1^- \neq 0$ by assumption. Thus, $g_1=0$.

\begin{case}
$f_0\in\B_0,g_1\in\B_1$ satisfy $f_0\B g_1=0$.
\end{case}  
Write as usual $f_0 = f_0^+  + f_0^- xy$ and $g_1 = g_1^+x  + g_1^- y$ where 
$f_0^{\pm},g_1^{\pm}\in\J$.
We have $f_0\B g_1x=0$; applying Case \ref{case1} to the even elements $f_0$ and $g_1x$, we see that either $f_0=0$ and we are done, or $f_0 \neq 0$. In the latter case,
$g_1x=0$, hence $g_1=g_1^+x$. 
The condition $f_0\B g_1=0$ implies in particular $(f_0^++f_0^-xy)g_1^+x=0$ in $\B_1$ hence $f_0^+g_1^++f_0^-\eta(g_1^+)z=0$ in $\J$. Consider $f_0yg_1^+x=0$; then
\begin{align*}
 0 &= (f_0^++f_0^-xy)yg_1^+x =(f_0^++f_0^-xy)(g_1^+y+\nabla (g_1^+)x)x\\
&=(f_0^++f_0^-xy)g_1^+yx=f_0^+g_1^+yx+f_0^-\eta(g_1^+)xy^2x\\
&=f_0^+g_1^+yx = f_0^+g_1^+z - f_0^+g_1^+xy
\end{align*}and deduce that $f_0^+g_1^+=0$ which implies that $f_0^-\eta(g_1^+) = 0$ from the previous
argument.
Therefore $g_1^+=0$, that is, $g_1=0$.

\begin{case}
$f_1\in\B_1,g_0\in\B_0$ satisfy $f_1\B g_0=0$. 
\end{case} 
Write as usual $f_1 = f_1^+ x + f_1^- y$ and $g_0 = g_0^+  + g_0^- xy$ where $f_1^{\pm},g_0^{\pm}\in\J$.
We have $xf_1\B g_0=0$; applying Case \ref{case1} to the even elements $xf_1$ and $g_0$, we see that either $g_0=0$ and we are done, or $g_0 \neq 0$. In the latter case,
 we have $xf_1=0$, i.e., $f_1=f_1^+x$. 
 We have the equalities $f_1^+x(g_0^++g_0^-xy)=0$ in $\B_1$ and  $f_1^+xy(g_0^++g_0^-xy)=0$ in $\B_0$. The first one leads to $f_1^+\eta(g_0^+)=0$ and the second one to 
$f_1^+\eta(g_0^+)+f_1^+\eta(g_0^-)z=0$. Thus  $f_1^+=0$ and  $f_1=0$.
\epf

We are now ready to localize the superalgebra
$\B$ at the set $S$,
extending our previous  Theorem \ref{QB0}.

\begin{theorem}\label{QB}
With the notations \emph{\ref{notaS}} we have the following properties:
\begin{enumerate}[leftmargin=*,label=\rm{(\roman*)}]
\item 
The set $\reg $ of homogeneous regular elements in $\B$ satisfies the left and right Ore conditions in $\B$,

\medbreak
\item  the ring of fractions $Q\coloneqq \B\reg ^{-1}$ is equipped with a structure of superalgebra whose even part is $Q_0=\B_0\reg _0^{-1}$ and we have $Q=\B\reg _0^{-1}$,

\medbreak
\item
the superalgebra $Q$ is a super-simple super-artinian ring.   
\end{enumerate}
\end{theorem}

\pf We have seen that $\B$ is $e$-faithful.
Now $\B_0$ is a (left and right) Goldie ring since it is Noetherian. 
Therefore assertions (i) and (ii) follow from Proposition \ref{B0sp} by \cite[Theorem 29]{BKM}. Then Theorem \ref{Bsp} implies assertion (iii) using \cite[Theorem 30]{BKM} or assertion (6) of \cite[Theorem 35]{BKM}.
\epf


\section{Super-automorphisms}

The goal of this section is to determine the group of super-automorphims of $\B$, that is the automorphisms of $\B$ preserving its superalgebra structure. 

\begin{lemma}\label{autB0lem1}Let $\sigma$ be an automorphism of the algebra $\B_0$. Then the commutative subalgebra $\ku z+\ku xy$ of $\B_0$ is stable by $\sigma$. More precisely there exist $\alpha,\beta\in\ku^\times$ such that either of the following possibilities hold:
\begin{align}\label{twotypes2}
&\left\{\begin{aligned} \sigma(z) &=\beta z +(\alpha - \beta) xy,
\\\sigma(xy) &=\alpha xy; \end{aligned}\right.
\\ \label{twotypes}
&
\left\{\begin{aligned}  
\sigma(z) &=\alpha z +(\beta - \alpha) xy,
\\\sigma(xy) &= \alpha z - \alpha xy. \end{aligned}\right.
\end{align}
\end{lemma}

\pf By \eqref{relxyz} we have $\sigma(xy)\sigma(z-xy)=0$.
We deduce from  \eqref{eq:zerodiv} that
there exist  $f_0^-$, $h_0^- \in \J$ such that
\begin{align}\label{twotypesbis-2}
\begin{aligned}
 &\text{either }& &{\rm (I)} &
 &\left\{\begin{aligned} \sigma(xy) &= f_0^-xy,\hfill\\
  \sigma(z-xy) &=h_0^- (z - xy), \\
 \sigma(z) &=h_0^- z+ (f_0^- -h_0^-) xy;
\end{aligned}\right.
\\ 
&\text{or} &&{\rm (II)} & 
& \left\{\begin{aligned} 
\sigma(xy) &= f_0^- (z -xy)\hfill\\
  \sigma(z-xy) &=h_0^-  xy, \\
 \sigma(z) &= f_0^-z + (h_0^- - f_0^-)xy.
\end{aligned}\right.
\end{aligned}
\end{align}

Since $\sigma^{-1}$ is also an automorphism, we may consider various cases.

\begin{stepp}\label{stepp1}
 If $\sigma$ is of type (I), then $\sigma^{-1}$  is of type (I).
\end{stepp}
 
 Let $f_0^-$, $h_0^- \in \J$ such that (I) holds for $\sigma$.
 Suppose that $\sigma^{-1}$  is of type (II).
 Then there exists $t_0^- \in \J$ such that
 $\sigma^{-1}(xy)= t_0^-(z - xy)$. 
 Now there exist $a_0^{\pm}\in\J$ such that
  $\sigma(t_ 0^-)=a_0^++a_0^-xy$. We compute
\begin{align*}
xy&=\sigma\sigma^{-1}(xy)=
\sigma(t_0^-)\sigma(z- xy)\\
&= (a_0^++a_0^-xy)h_0^-(z-xy)\\
&= a_0^+h_0^-(z-xy)+a_0^-\eta(h_0^-)xy(z-xy)\\
&= a_0^+h_0^-(z-xy)
\end{align*}
which is impossible in $\B_0=\J\oplus\J xy$.

\begin{stepp}\label{stepp2}
If $\sigma$ is of type (I), then 
there exist $\alpha,\beta\in\ku^\times$ such that
\eqref{twotypes2} holds.
\end{stepp} 

 By Claim \ref{stepp1}, $\sigma^{-1}$ is necessarily of type (I). There exist $t_0^-, a_0^{\pm}\in \J$ such that
 $\sigma^{-1}(xy)=t_0^-xy$  and $\sigma(t_ 0^-)=a_0^++a_0^-xy$. Hence
\begin{align*}
xy&=\sigma\sigma^{-1}(xy)=\sigma(t_0^-)\sigma(xy)=(a_0^++a_0^-xy)f_0^-xy\\
&=\bigl(a_0^+f_0^-+a_0^-\eta(f_0^-)z\bigr)xy=(a_0^++a_0^-z)f_0^-xy.
\end{align*}
This equality in $\B_0$ implies $(a_0^++a_0^-z)f_0^-=1$ in $\J$ then $f_0^-=\alpha$ for some $\alpha\in\ku^\times$. Hence $\sigma(xy)=\alpha xy$, which implies $\sigma^{-1}(xy)=\alpha^{-1}xy$.

\medbreak Also, there exist $s_0^-, b_0^{\pm}\in \J$ such that $\sigma^{-1}(z - xy) = s_0^-(z - xy)$  and 
$\sigma(s_ 0^-)=b_0^+ + b_0^-xy$. Hence
\begin{align*}
z - xy&= \sigma\sigma^{-1}(z - xy) = \sigma(s_0^-)\sigma(z - xy) = (b_0^+ + b_0^-xy)h_0^-(z - xy)\\
&= b_0^+ h_0^-(z - xy) + b_0^- \eta(h_0^-) xy(z - xy) = b_0^+ h_0^-(z - xy)
\end{align*}
 This equality in $\B_0$ implies $b_0^+ h_0^-=1$. Thus
 there exists $\beta\in\ku^\times$ such that  $\sigma(z - xy)=\beta (z-xy)$. Thus we are in the situation \eqref{twotypes2}.

\begin{stepp}\label{stepp3}
There exists an involutive automorphism $\omega$ 
of $\B_0$, which is of type (II),  defined~by
\begin{align}\label{eq:omega}
\omega(z) &= z, & \omega(y^2)&=y^2,&  
\omega(xy) &= z-xy. 
\end{align}
\end{stepp}

The fact that the linear isomorphism $\omega$
is multiplicative follows from the compatibility of the assignements \eqref{eq:omega} with the identities in 
Remark \ref{rem:algebra-B_0}. 
Clearly,   $\omega^2=\id$; taking
$g_0^-=0$ and $f_0^-=-1$ in \eqref{twotypes2} we see that 
$\omega$  is of type (II). 

\begin{stepp}\label{stepp4}
If $\sigma$ is of type (II), then 
there exist $\alpha,\beta\in\ku^\times$ such that
\eqref{twotypes} holds.
\end{stepp}

Then $\omega\sigma$ is of type (I),
hence   there exist $\alpha,\beta\in\ku^\times$ such that 
\begin{align*}
\omega\sigma (z) &=  \beta z+(\alpha-\beta)xy, &
\omega\sigma (xy) &= \alpha xy
\end{align*}
by Claim \ref{stepp2}.
Therefore
\begin{align*}
\sigma(z) &= \omega(\beta z+(\alpha-\beta)xy)=\beta z+(\alpha-\beta)(z-xy)=\alpha z + (\beta - \alpha)xy,
\\
 \sigma(xy) &= \omega(\alpha xy)= \alpha z-\alpha xy.
\end{align*}
The Lemma is proved. \epf

\begin{lemma}\label{lema:autom}
\emph{(i)} Given $\mu\in\ku^\times$ and $s(z)\in\ku[z]$,
there exists a unique super-automorphism $\sigma_{s,\mu}$ 
of $\B$  such that
\begin{align}\label{sigmax&y}
\sigma_{s,\mu}(x) &=\mu x&&\text{and}& 
\sigma_{s,\mu}(y) &=\mu y+s(z)x.
\end{align}

\emph{(ii)} If also  $\mu'\in\ku^\times$ and 
$s'(z)\in\ku[z]$, then
\begin{align}\label{composigma}
\sigma_{s, \mu}\sigma_{s', \mu'} &=\sigma_{s'',\mu\mu'}& 
&\text{where}&  s''(z) &= \mu' s(z)+\mu s'(\mu^2z).
\end{align}
\end{lemma}
\pf (i) follows from the defining relations
\eqref{relB};  the proof of (ii) is direct.
\epf

\begin{theorem}\label{AutB}Let $\sigma$ be a super-automorphism of  $\B$.
\begin{enumerate}[leftmargin=*,label=\rm{(\roman*)}]
\item
There exist $\mu\in~\ku^\times$ and $s(z)\in\ku[z]$ such that $\sigma = \sigma_{s,\mu}$.

\medbreak
\item\label{item:N} The group $\Auts \B$
of superalgebra automorphisms of $\B$ is isomorphic to 
the semidirect product $N \rtimes K$, 
where $K \simeq \Bbbk^\times$ acts on $N \simeq \Bbbk[z]$ by $\mu \cdot s(z) = s(\mu^2z)$;
concretely, the isomorphism $N \rtimes K \to \Auts \B$ is given by $(s, \mu) 
\mapsto \sigma_{\mu s, \mu}$, $s \in \Bbbk[z]$, $\mu \in \Bbbk^{\times}$. Thus
\begin{align*}
N&\to \{\sigma_{s,1}: s\in\Bbbk[z]\},&
K &\to \{\sigma_{0,\mu}\,;\,\mu\in\Bbbk^\times\}.
\end{align*}

\medbreak
\item The restriction of $\sigma$ to $\B_0$ is given by
\begin{equation}\label{sigmazxyy2}
\sigma(z)=\mu^2z,\quad \sigma(y^2)=\mu^2 y^2+\mu zs(z),\quad\sigma(xy)=\mu^2xy.
\end{equation}
\end{enumerate}\end{theorem}

\pf  There exist $f_1^{\pm}$,  $g_1^{\pm}\in\J$ such that $\sigma(x)=f_1^+x+f_1^-y$ and 
$\sigma(y)=g_1^+x+g_1^-y$. 
Since $\sigma(x)^2=0$ it follows from 
Proposition \ref{oddZD} that either
$f_1^-=0$ or  $f_1^-y^2=-f_1^+z$ with $f_1^-\not=0$.

\begin{stepo}
If $f_1^-=0$, then $\sigma = \sigma_{s,\mu}$
for unique $\mu\in~\ku^\times$ and $s(z)\in\ku[z]$.
\end{stepo}
Suppose that $f_1^-=0$. We compute using relation \eqref{twoodd} and adding:
\begin{align*} 
\sigma(xy)&=\sigma(x)\sigma(y) =f_1^+\eta(g_1^-)xy
\\ \notag
\sigma(yx)&=\sigma(y)\sigma(x)
=g_1^-f_1^+z+g_1^-(-f_1^++\nabla (f_1^+))xy,
\\ \notag
\sigma(z) &= g_1^-f_1^+z+\bigl(f_1^+\eta(g_1^-)-g_1^-f_1^++g_1^-\nabla (f_1^+)\bigr)xy.
\end{align*}
By Lemma \ref{autB0lem1} and inspection we see
that there exists $\alpha, \beta \in \ku^\times$ such that  \eqref{twotypes2} holds.  Then there exist $\lambda,\mu\in\ku^\times$ such that 
$f_1^+=\lambda$ and $g_1^-=\mu$.
Clearly,  $\lambda\mu=\beta$ and $\alpha-\beta=0$. We have 
\begin{align*}
\sigma(x) &=\lambda x, & \sigma(y) &=g_1^+x+\mu y,&
\sigma(xy)&=\lambda\mu xy, &\sigma(z)&=\lambda\mu z.
\end{align*}

Applying $\sigma$ to relation \eqref{relB} we  obtain
\begin{align*}
0 &= \sigma(y)\sigma(z) -\sigma(z)\sigma(y)- \sigma(x)\sigma(z)
\\
&= (g_1^+x+\mu y)\lambda\mu z-\lambda\mu z(g_1^+x+\mu y)-\lambda x\lambda\mu z
\\
&= (g_1^+z-zg_1^+)x+\mu(yz-zy)-\lambda xz
\\
&= (g_1^+z- \eta(g_1^+)z)x + (\mu -\lambda) zx,
\end{align*}

Therefore $g_1^+-\eta(g_1^+)+(\lambda-\mu)=0$. 
 
 By formula \eqref{defetabis} we necessarily have  $\lambda=\mu$. Then $\eta(g_1^+)=g_1^+$. 
To sum up $\sigma(x)=\mu x$ and  $\sigma(y)= sx + \mu y$,  where $s=g_1^+\in\J^\eta$.

When $\ch\ku=0$ we have $s\in\ku[z]$ by \eqref{eq:defeta} and we obtain  \eqref{sigmax&y}.

Assume that $\ch\ku=p>2$. We have $s\in\ku[z,y^{2p}]$  by \eqref{eq:defeta}. Using $sx=xs$ and $x^2=0$ we compute:
\begin{align*}
\sigma(y^2) & = (\mu y+sx)(\mu y + sx)
\\
&= \mu^2y^2+\mu sxy+\mu(sy+\nabla(s)x)x
\\
& =\mu^2y^2+\mu sz=\mu^2y^2+\mu zs.
\end{align*}
The restriction of $\sigma$ to $Z(\B)$ is an automorphism of $Z(\B)$. We have $Z(\B)=\ku[z^p,y^{2p}]$ by Proposition \ref{centerB} and $\sigma(z^p)=\lambda^p\mu^pz^p$. 
Hence there exist $\alpha\in\ku^\times$ and $t\in\ku[z^p]$ such that $\sigma(y^{2p})=\alpha y^{2p}+t$. The equality $(\mu^2y^2+\mu zs)^p=\alpha y^{2p}+t$ implies $s\in\ku[z]$. We obtain  \eqref{sigmax&y}.

\begin{stepo}
There is no automorphism satisfying $f_1^-\not=0$.
\end{stepo}
Suppose that $f_1^-\not=0$; hence $f_1^-y^2=-f_1^+z$. We compute using \eqref{twoodd}:
\begin{align*}
\sigma(xy)&=f_1^-(g_1^+z+g_1^-y^2)+\Bigl(f_1^+\eta(g_1^-)+f_1^-\bigl(\nabla (g_1^-)-g_1^+\bigr)\Bigr)xy
\\
\sigma(yx)&=g_1^-(f_1^+z+f_1^-y^2)+\Bigl(g_1^
+\eta(f_1^-)+g_1^-\bigl(\nabla (f_1^-)-f_1^+\bigr)\Bigr)xy
\\
&= \Bigl(g_1^
+\eta(f_1^-)+g_1^-\bigl(\nabla (f_1^-)-f_1^+\bigr)\Bigr)xy
\\
\sigma(z) &=\sigma(xy)+\sigma(yx) \in f_1^-(g_1^+z+g_1^-y^2) + \J xy.
\end{align*}
By Lemma \ref{autB0lem1}, there exists $\gamma\in\ku^\times$ such that $f_1^-(g_1^+z+g_1^-y^2)=\gamma z$. 
Since $\J$ is an Ore extension, by Proposition \ref{B+},
we see  that $f_1^-\in\ku[z]$, $g_1^+z+g_1^-y^2\in\ku[z]$ 
and $\deg_z(f_1^-)+\deg_z(g_1^+z+g_1^-y^2)=1$.

\medbreak
If $f_1^- = \lambda\in\ku^{\times}$, then 
$f_1^+=a_0(z)+a_1(z)y^2$ where $a_0(z),a_1(z)\in\ku[z]$,
since $f_1^-y^2=-f_1^+z$.
Comparing coefficients in $\J=\ku[z][y^2\,;\, d^2]$,  $\lambda=-a_1(z)z$, which is impossible.

\medbreak
Assume then that $\deg_z(f_ 1^-)=1$ and $g_1^+z+g_1^-y^2\coloneqq \lambda'\in\ku^\times$.
The coefficient of~$y^0$ in the expansion in $\J$ 
of the left hand is the coefficient of $y^0$ 
in the expansion of $g_1^+z$, 
which belongs to the ideal $z\ku[z]$ by \eqref{commJ}. This could not be equal  to $\lambda'$,  again a contradiction. Thus the case $f_1^-\not=0$ is impossible. 

\medbreak 
Thus (i) is proved; then (ii) follows  from Lemma \ref{lema:autom} and (iii)  from~(i).
\epf

\begin{remark}\label{compAutBJ}
Recall from \cite[Proposition 3.6]{AlD}, \cite[Theorem 2.7, Theorem 2.10]{Sh2007}, \cite[Theorem 8.2]{BLO}, that the automorphisms of $\J$ have the form
\begin{align*}
\tau_{p,\alpha}:\, z&\mapsto \alpha z, &y^2 &\mapsto\alpha y^2+p(z), & \alpha &\in\ku^\times,\quad p(z) \in\ku[z],
\end{align*} with composition
\begin{equation*}
\tau_{p,\alpha}\tau_{p',\alpha'}=\tau_{p'',\alpha\alpha'}\qquad \text{where \ } p''(z)=\alpha' p(z)+p'(\alpha z).
\end{equation*}
By \eqref{sigmazxyy2}  $\J$ is stable by any super-automorphism $\sigma_{s,\mu}$ of $\B$, its restriction being $\tau_{\mu zs,\mu^2}$; i.e., we have a group homomorphism $\Res: \Auts \B \to \Aut \J$ and
\begin{align*}
\ker \Res &= \{\sigma_{0, \pm 1}\}, &
\imm \Res &= \{\tau_{p,\alpha}: p(z) \in z\ku[z], \quad
\alpha\in (\ku^{\times})^2\}.
\end{align*}
\end{remark}


\section{On the derivations of $\B$} 

As is well-known, the first Hochschild homology group $\mathrm{HH}^1(\B)$
is isomorphic to the module of outer derivations of $\B$, i.e.,
$\mathrm{HH}^1(\B) \simeq \Der \B / \Derinn B$. The Lie algebra
$\mathrm{HH}^1(\B)$ was computed explicitly in \cite{RS} when $\ch\ku=0$:

\begin{theorem}\label{derB}\cite[Theorem 3.12, Propositions 6.1 and 6.2, Theorem 6.3]{RS}
Assume that $\ch \ku =0$. The first Hochschild homology group $\mathrm{HH}^1(\B)$ has a basis 
$\{\mathfrak c\} \cup \{\mathfrak s_n: n \in \N_0\}$ where the Lie bracket satisfies
\begin{align*}
[\mathfrak c, \mathfrak s_n] &= 0, &[\mathfrak s_m, \mathfrak s_n] &= 2(n - m)\mathfrak s_{n+m}, &m, n &\in \N_0.
\end{align*}
Thus, $\mathrm{HH}^1(\B)$ is isomorphic to a subalgebra of the Virasoro algebra. \qed
\end{theorem}
 
 Recall the subgroup $N$ of $\Auts \B$ introduced in Theorem \ref{AutB} \ref{item:N}.
\begin{prop} \label{prop:locnilp}
Assume that $\ch \ku =0$. Any automorphism of $\B$ belonging to 
 $N$ is the exponential of a locally nilpotent derivation of $\B$.
\end{prop}
\pf
We introduce the $\ku$-derivation $\mathfrak c$ of $\B$   by $\mathfrak c(x)=0$ and $\mathfrak c(y)=x$; then $\mathfrak c(xy)=\mathfrak c(z)=0$; thus $\ku[z]\subseteq\ker \mathfrak c$.  By elementary calculations
using \eqref{BOre} we check that $\deg_y(\mathfrak c(b))\leq\deg_y(b)-1$ for any $b\in\B$. Hence $\mathfrak c$ is locally nilpotent. 
Let $s=\sum_{i\geq 0}\alpha_iz^i \in\ku[z]$, $\alpha_i\in\ku$.
Then $s\mathfrak c$ is a locally nilpotent derivation of $\B$. 
Notice that the image of $\mathfrak c$ in $\mathrm{HH}^1(\B)$ 
is the element with the same name.
Now, the derivation $s\mathfrak c$ decomposes in the basis mentioned in Theorem \ref{derB} as
\begin{equation*}s\mathfrak c=s(0)\mathfrak c+\rm{ad}_{\overline s}
\end{equation*}
where $\overline{s}=\sum_{i\geq 1}\frac{\alpha_i}iz^i$.
We have $\sigma_{s,1}=\exp(s\mathfrak c)$.\epf
The derivation $\mathfrak c$ satisfies $\mathfrak c(y^2)=z$ and denoting by $\widetilde{\mathfrak c}$ its restriction to $\J$ we have  $\exp(s\widetilde{\mathfrak c})=\tau_{zs,1}$ with the notations of Remark \ref{compAutBJ}.


\section{The super Jordan plane as a braided Hopf algebra}
\label{section:braided-Hopf}

Recall that $p =\ch \ku \neq 2$. 
Consider the braided vector space $(V,c)$ 
where $V=\ku x\oplus\ku y$ has dimension 2
and   the braiding $c:V\otimes V\to V\otimes V$ is given~by
\begin{align*}
c(x\otimes x) &=- x\otimes x, &c(y\otimes x)&=- x\otimes y,
\\ c(x\otimes y)&=(-y+x)\otimes x & 
c(y\otimes y) &= (- y+x)\otimes y. 
\end{align*}
Fix a generator $g$ of the group $\Z$.
Then $(V,c)$ can be realized as a Yetter-Drinfeld module
over $\Bbbk \Z$ by 
\begin{align}\label{eq:realisation}
\delta(v) &= g \otimes v,& v&\in V, &g\cdot x &= - x,
& g\cdot y&= -y + x.
\end{align}
Then $\B$ is a Hopf algebra in $\yd{\ku \Z}$
 by declaring that $x$ and $y$ are primitives, see \cite[Section~3.3]{AAH2}. 
 Actually, we know:

\begin{itemize}[leftmargin=*]\renewcommand{\labelitemi}{$\circ$}

\medbreak
\item If $p=0$, then $\B$ is isomorphic to
the Nichols algebra $\toba(V)$,  \cite[Section~3.3]{AAH2}.

\medbreak 
\item If $p > 2$, then  
$\toba(V) \simeq \B/ \langle y^{2p}, z^p\rangle$, see \cite[Section 3.2]{AAH3}. Indeed we have:
\end{itemize}

\begin{lemma}
If $p > 2$, then there exists an exact sequence of braided Hopf algebras
\begin{align}\label{eq:exact-seq-nichols}
\ku  [y^{2p}, z^p] \hookrightarrow \B \overset{\varpi}{\twoheadrightarrow} \toba(V).
\end{align}
\end{lemma}

\pf First we claim that the algebra $\ku  [y^{2p}, z^p]$, which 
 is $Z(\B)$ as we know, is a Hopf subalgebra of $\B$. 
Indeed, $\deg y^{2p} = \deg z^p = 2p$, i.e., $y^{2p}$ and $z^p$ are elements of the smallest degree in the Hopf ideal $\langle y^{2p}, z^p\rangle$; hence they are primitive (alternatively one may use
 the formulas \cite[(3.12)]{APP} to show that $y^{2p}$ and $z^p$ are primitive). 
 Notice that $ g^n\cdot y = (-1)^n (y - nx)$ for all $n \in \N$, hence
$ g^{2p}\cdot y^j = y^j$ for all $j \in \N$. Also,
\begin{align*}
g\cdot z   &= g\cdot (xy + yx) = -x(-y +x) + (-y +x)(-x) = z. 
\end{align*}
Therefore the restriction of the braiding to
$\ku  [y^{2p}, z^p]$ is  the usual transposition;
together with the claim that $y^{2p}$ and $z^p$ are primitive,
this implies that $\ku  [y^{2p}, z^p]$ is a braided Hopf subalgebra
of $\B$. Finally, we have to check that 
 $\ker \varpi = \ku  [y^{2p}, z^p]^+ \B$ and 
$\B^{\co \varpi} =  \ku  [y^{2p}, z^p]$
(where $\ku  [y^{2p}, z^p]^+$ is the kernel of the counit of $\ku  [y^{2p}, z^p]$).
We omit  these verifications that are standard.
\epf

\begin{remark}
Observe that $A = \ku[x]/(x^2)$ is a Hopf superalgebra by declaring that 
$x$ is primitive. Thus the inclusion $\iota: A \to \B$ is morphism of braided Hopf algebras. Now there is an algebra map $\omega :\B \to A$ given by 
$\omega(x) = x$ and $\omega(y) = 0$.  
However it is not a morphism of braided Hopf algebras because it does not preserve
 the braiding; even if $x$ and $y$ are primitives in $\B$, the map $\omega$
 is not a morphism of coalgebras. 
 \end{remark}

\begin{prop}\label{prop:exseq}
The braided Hopf algebra $\B$ fits into an extension
\begin{align}\label{eq:exseq}
A[z] \hookrightarrow \B \overset{\pi}{\twoheadrightarrow} L
\end{align}
 of braided Hopf algebras,
where $A[z] = \ku  [ x, z]$ as defined above; 
$L = \ku[\eta]$, a polynomial ring with
braiding and comultiplication determined by 
$c(\eta \otimes \eta) = -\eta \otimes \eta$ 
and $\Delta(\eta) = \eta \otimes 1 + 1 \otimes \eta$ respectively;
and $\pi$ is determined by
\begin{align*}
\pi(x) &= 0 &&\text{and}& \pi(y) &= \eta.
\end{align*}
\end{prop}

\pf The subalgebra $A[z]$ is commutative by \eqref{relBbis};
it is a Yetter-Drinfeld submodule of $\B$ with respect to 
the realization \eqref{eq:realisation} because 
$\delta (z) = g^2 \otimes z$ and $g \cdot z = z$.
Hence the inclusion preserves the braiding. Therefore $A[z]$
is also a subcoalgebra since $x$ is primitive and $\Delta (z) = z \otimes 1 + 1 \otimes z + x \otimes x$. It follows at once from \eqref{relB} that $A[z]$ is normal.
Now $L$ is a Hopf algebra in $\yd{\ku \Z}$ by declaring 
$\delta(\eta) = g \otimes \eta$ and $g\cdot \eta = - \eta$. 
Clearly, $\pi$ is an algebra map and a morphism 
of Yetter-Drinfeld modules, hence it preserves the braiding and a fortiori 
it is a coalgebra map. Hence \eqref{eq:exseq} is an exact sequence of Hopf algebras.
\epf
Keep the notation above.

\begin{prop}Assume that $p>2$.
The  Nichols algebra $\toba(V)$ fits into an extension of Hopf algebras
\begin{align}\label{eq:exseq-p}
A[z] / \langle z^p\rangle\hookrightarrow \toba(V) \overset{\pi}{\twoheadrightarrow}  L/\langle \eta^{2p}\rangle.
\end{align}
\end{prop}
\pf
This can be derived from Proposition \ref{prop:exseq}, or proved analogously.
\epf


\section{Hopf algebra automorphisms}
Let $\Aut_{\text{Hopf}}(\B)$ be the group of (braided)
Hopf algebra automorphisms of $\B$ and 
$G \coloneqq \{\sigma_{s,\mu} \in \Aut \B: \mu \in \ku^{\times}, \,  s\in\ku\}$, a subgroup of $\Auts (\B)$.

\begin{prop} \label{prop:hopf-automorphisms}
We have $\Auth(\B) \cap \Auts (\B) = G$. Furthermore, 
if $\ch \ku = 0$, then $\Auth(\B) = G$.
\end{prop}

\pf  (i) If $\mu \in \ku^{\times}$ and $s\in \ku$, then
 $\sigma_{s,\mu}(x), \sigma_{s,\mu}(y) \in V = \ku x \oplus \ku y \subseteq \Pc (\B)$.
A direct computation shows that
$\sigma_{s,\mu} \otimes \sigma_{s,\mu}$ 
commutes with the braiding $c$. Therefore, $\sigma_{s,\mu}$
preserves the comultiplication and is thus a  Hopf algebra map.
This shows that $\Auth(\B) \cap \Auts (\B) \supset G$.

\medbreak
(ii) Conversely, let $\sigma \in \Auth(\B) \cap \Auts (\B)$.
By Theorem \ref{AutB},  $\sigma = \sigma_{s,\mu}$
for unique $\mu\in\ku^\times$ and $s(z)\in\ku[z]$. 
Then $\sigma_{s,\mu}(y) =\mu y+s(z)x \in \Pc(\B)$,
hence $u \coloneqq s(z)x \in \Pc(\B)$. 
Assume that $\ch \ku > 2$ and let $\varpi$ be as in \eqref{eq:exact-seq-nichols}. 
Then
$\varpi(u)  \in \Pc(\toba(V)) = V$, i.e., there exist
$\alpha, \beta\in\ku$ such that
$w \coloneqq u - \alpha x -\beta y \in \Pc(\B) \cap \ker \varpi$.
Hence
\begin{align*}
(\id \otimes \varpi )\Delta (w) = w \otimes 1 \implies
w \in \B^{\co \varpi}  = \ku  [x, z].
\end{align*}
That is,  
$s(z)x - \alpha x -\beta y \in  \ku  [x, z]$. This implies that $s(z) \in \ku$.

\medbreak
(iii) Assume that $\ch \ku = 0$; then  \cite[Theorem 1.4]{AD-irma}
applies. For the sake of completeness we repeat the argument.
Let $\sigma \in \Aut_{\text{Hopf}} (\B)$. Then $\sigma$, being a coalgebra map, preserves  $\Pc(\B)$ which  is $V$  since $\B \simeq \toba(V)$; 
particularly, 
$\sigma(x), \sigma(y) \in V \subset \B_1$. Thus $\sigma$ is a super-automorphism, since $\B_0$, respectively
$\B_1$,  is spanned by the monomials in $x$ and $y$ of
even, respectively odd, length.
Hence the second statement follows from the first.
\epf


\subsection*{Acknowledgements}
We thank Andrea Solotar for interesting exchanges on the question of the global dimension.

 \bibliographystyle{plain}
 \bibliography{superJordan}

\begin{thebibliography}{10}

\bibitem{AlD}
Jacques Alev and Fran\c{c}ois Dumas.
\newblock Invariants du corps de {W}eyl sous l'action de groupes finis.
\newblock {\em Comm. Algebra}, 25(5):1655--1672, 1997.

\bibitem{AAH1}
Nicol\'as Andruskiewitsch, Iv\'an Angiono, and Istv\'an Heckenberger.
\newblock Liftings of {J}ordan and super {J}ordan planes.
\newblock {\em Proc. Edinb. Math. Soc. (2)}, 61(3):661--672, 2018.
\newblock Corrigendum. \textit{Proc. Edinb. Math. Soc. (2)}, 65(3):577--586,
  2022.

\bibitem{AAH2}
Nicol\'as Andruskiewitsch, Iv\'an Angiono, and Istv\'an Heckenberger.
\newblock On finite {GK}-dimensional {N}ichols algebras over abelian groups.
\newblock {\em Mem. Amer. Math. Soc.}, 271(1329):ix+125, 2021.

\bibitem{AAH3}
Nicol\'as Andruskiewitsch, Iv\'an Angiono, and Istv\'an Heckenberger.
\newblock Examples of finite-dimensional pointed {H}opf algebras in positive
  characteristic.
\newblock In {\em Representation theory, mathematical physics, and integrable
  systems}, volume 340 of {\em Progr. Math.}, pages 1--38.
  Birkh\"auser/Springer, Cham, [2021] \copyright 2021.

\bibitem{ABDF17}
Nicol\'as Andruskiewitsch, Dirceu Bagio, Saradia Della~Flora, and Daiana
  Fl\^ores.
\newblock Representations of the super {J}ordan plane.
\newblock {\em S\~ao Paulo J. Math. Sci.}, 11(2):312--325, 2017.

\bibitem{ABDF}
Nicol\'as Andruskiewitsch, Dirceu Bagio, Saradia Della~Flora, and Daiana
  Fl\^ores.
\newblock On the bosonization of the super {J}ordan plane.
\newblock {\em S\~ao Paulo J. Math. Sci.}, 13(1):1--26, 2019.

\bibitem{AD-irma}
Nicol\'as Andruskiewitsch and Fran\c{c}ois Dumas.
\newblock On the automorphisms of {$U_q^+(g)$}.
\newblock In {\em Quantum groups}, volume~12 of {\em IRMA Lect. Math. Theor.
  Phys.}, pages 107--133. Eur. Math. Soc., Z\"urich, 2008.

\bibitem{APP}
Nicol\'as Andruskiewitsch and H\'ector~Mart\'in Pe\~na Pollastri.
\newblock On the double of the (restricted) super {J}ordan plane.
\newblock {\em New York J. Math.}, 28:1596--1622, 2022.

\bibitem{ArS}
Michael Artin and William~F. Schelter.
\newblock Graded algebras of global dimension {$3$}.
\newblock {\em Adv. in Math.}, 66(2):171--216, 1987.

\bibitem{BKM}
Irina~N. Balaba, Andrei~L. Kanunnikov, and Aleksandr~V. Mikhalev.
\newblock Quotient rings of graded associative rings. {I}.
\newblock {\em J. Math. Sci., New York}, 186(4):531--577, 2012.
\newblock {T}ranslated from \textit{Fundam. Prikl. Mat.}, {17} (2):3--74,
  2011/12.

\bibitem{BLO}
Georgia Benkart, Samuel~A. Lopes, and Matthew Ondrus.
\newblock A parametric family of subalgebras of the {W}eyl algebra {I}.
  {S}tructure and automorphisms.
\newblock {\em Trans. Amer. Math. Soc.}, 367(3):1993--2021, 2015.

\bibitem{CS}
Javier C\'oppola and Andrea Solotar.
\newblock Graded braided commutativity in {H}ochschild cohomology.
\newblock {\em Theory Appl. Categ.}, 41:Paper No. 46, 1596--1643, 2024.

\bibitem{DMMZ}
Evgeny~E. Demidov, Yuri \~I. Manin, Evgeny~E. Mukhin, and Dimitri~V.
  Zhdanovich.
\newblock Nonstandard quantum deformations of {${\rm GL}(n)$} and constant
  solutions of the {Y}ang-{B}axter equation.
\newblock {\em Progr. Theoret. Phys. Suppl.}, 102:203--218, 1990.
\newblock Common trends in mathematics and quantum field theories (Kyoto,
  1990).

\bibitem{Dix}
Jacques Dixmier.
\newblock {\em Alg\`ebres enveloppantes}, volume Fasc. XXXVII of {\em Cahiers
  Scientifiques}.
\newblock Gauthier-Villars \'Editeur, Paris-Brussels-Montreal, Que., 1974.

\bibitem{GW}
Kenneth~R. Goodearl and Robert~B. Warfield, Jr.
\newblock {\em An introduction to noncommutative {N}oetherian rings}, volume~61
  of {\em London Mathematical Society Student Texts}.
\newblock Cambridge University Press, Cambridge, second edition, 2004.

\bibitem{G}
Dimitri~I. Gurevich.
\newblock The {Y}ang-{B}axter equation and the generalization of formal {L}ie
  theory.
\newblock {\em Dokl. Akad. Nauk SSSR}, 288(4):797--801, 1986.

\bibitem{Iyudu}
Natalia~K. Iyudu.
\newblock Representation spaces of the {J}ordan plane.
\newblock {\em Comm. Algebra}, 42(8):3507--3540, 2014.

\bibitem{LM}
Andrey~Yu. Lazarev and Mikhail~V. Movshev.
\newblock Quantization of some {L}ie groups and algebras.
\newblock {\em Uspekhi Mat. Nauk}, 46(6(282)):215--216, 1991.

\bibitem{mcconnell-robson}
John~C. McConnell and J.~Chris Robson.
\newblock {\em Noncommutative {Noetherian} rings. {With} the cooperation of
  {L}. {W}. {Small}.}, volume~30 of {\em Grad. Stud. Math.}
\newblock Providence, RI: American Mathematical Society (AMS), reprint\-ed with
  corrections from the 1987 original edition, 2001.

\bibitem{NvO1}
Constantin N\u{a}st\u{a}sescu and Freddy Van~Oystaeyen.
\newblock {\em Methods of graded rings}, volume 1836 of {\em Lecture Notes in
  Mathematics}.
\newblock Springer-Verlag, Berlin, 2004.

\bibitem{ohn}
Christian Ohn.
\newblock A {$*$}-product on {${\rm SL}(2)$} and the corresponding nonstandard
  quantum-{$U(sl(2))$}.
\newblock {\em Lett. Math. Phys.}, 25(2):85--88, 1992.

\bibitem{RS}
Sebasti\'an Reca and Andrea Solotar.
\newblock Homological invariants relating the super {J}ordan plane to the
  {V}irasoro algebra.
\newblock {\em J. Algebra}, 507:120--185, 2018.

\bibitem{Shirikov0}
Evgenij~N. Shirikov.
\newblock Two-generated graded algebras.
\newblock {\em Algebra Discrete Math.}, 2005(3):60--84, 2005.

\bibitem{Sh2007}
Evgenij~N. Shirikov.
\newblock The {J}ordanian plane.
\newblock {\em Fundam. Prikl. Mat.}, 13(2):217--230, 2007.

\bibitem{zak}
Stanisław Zakrzewski.
\newblock A {H}opf star-algebra of polynomials on the quantum {${\rm SL}(2,{\bf
  R})$} for a ``unitary'' {$R$}-matrix.
\newblock {\em Lett. Math. Phys.}, 22(4):287--289, 1991.

\end{thebibliography}
\end{document}